\documentclass[a4paper,12pt]{article}
\usepackage{amssymb,amsmath,enumerate,empheq,amsthm,amsfonts}
\usepackage[bottom]{footmisc}
\usepackage{dsfont}
\usepackage[dvips]{graphicx}

\usepackage{verbatim}
\usepackage{alltt}

\usepackage{color}

\usepackage[top=60pt, bottom=60pt, left=60pt, right=60pt]{geometry}

\usepackage{tikz}
\usetikzlibrary{matrix,arrows}

\usepackage{pifont} 

\usepackage{alltt} 
\usepackage{verbatim}

\usepackage{makeidx} 

\newtheorem{thm}{Theorem}[subsection]

\theoremstyle{definition}
\newtheorem{defn}[thm]{Definition}

\swapnumbers

\renewenvironment{abstract}
               {\list{}{\rightmargin\leftmargin}%
                \item[\hspace{10mm}\textbf{Abstract.}]\relax}
               {\endlist}

\allowdisplaybreaks

\makeindex
\begin{document}
\title{RANK OF THE DERIVATIVE OF THE PROJECTION TO SYMMETRIZED POLYDISC}   
\author{TRAN DUC ANH\footnote{Received ........ Revised ........ Accepted  ........\newline Contact Tran Duc Anh,  e-mail address: ducanh@hnue.edu.vn.}}

\date{\emph{Faculty of Mathematics, Hanoi National University of Education}\\ \emph{136 Xuan Thuy Str, Cau Giay, Hanoi}}

\maketitle

\abstract We prove that the rank of the derivative of the projection from spectral unit ball to symmetrized polydisc is equal to the degree of the minimal polynomial of the matrix at which we take derivative.\endabstract

\section{Introduction} In this note, we are interested in a special mapping in the spectral Nevanlinna-Pick problem which is called the symmetrization mapping. Firstly we present some notations. 

Denote by $\mathbb{C}^{n,n}$ the set of complex square matrices of size $n$ where $n$ is a positive integer. For each matrix $M\in \mathbb{C}^{n,n},$ the characteristic polynomial of $M$ is defined to be $$P_M(t) = \det (tI-M) = \sum_{j=0}^n (-1)^j\sigma_j(M)t^{n-j}$$ where $\sigma_0(M) = 1$ by convention, $\sigma_j(M)$ are the coefficients of the characteristic polynomial $\det (tI-M)$ and $I$ is the usual unit matrix. The $\sigma_j(M)'s$ are in fact the j-th elementary symmetric functions of the eigenvalues of $M.$

Put $\pi(M) = (\sigma_1(M), \sigma_2(M),\ldots, \sigma_n(M))$ so we get a mapping $\pi\colon \mathbb{C}^{n,n}\to \mathbb{C}^n$ which is called the \textit{symmetrization mapping}. 

Next we put $$\Omega_n = \{M\in \mathbb{C}^{n,n}~:~ r(M)<1\}$$ where $r(M)$ is the spectral radius of $M.$ We call $\Omega_n$ the \textit{spectral ball} and its image $$\mathbb{G}_n = \pi(\Omega_n)$$ the \textit{symmetrized polydisc}. This object is firstly introduced by J. Agler and N. Young \cite{Agler} to study spectral Nevanlinna-Pick problem, i.e., given $\alpha_1, \alpha_2, \ldots, \alpha_m\in \mathbb{D} = \{z\in \mathbb{C}~:~ |z|<1\}$ and $A_1, A_2,\ldots, A_m \in \Omega_n,$ find conditions such that there exists a holomorphic function $\Phi\colon \mathbb{D}\to \Omega_n$ with $\Phi(\alpha_j) = A_j$ for $1\leq j\leq m.$

Unfortunately, this problem is extremely difficult since $\Omega_n$ has an unfriendly geometry in comparison to the homogeneity of the disc or the unit ball in the classical Nevanlinna-Pick interpolation problem. So J. Agler and N. Young propose a new idea that instead of studying $\Omega_n,$ they project $\Omega_n$ onto $\mathbb{G}_n$ by symmetrization map and this object is a bounded taut domain, so theoretically it is easier to deal with. Since then, a lot of papers appear to study the symmetrized polydisc and give interesting information. 

In the inverse direction, we are interested in the lifting problem, i.e., if we can realize interpolation in $\mathbb{G}_n,$ can we lift it to $\Omega_n?$ A partial question is solved in \cite{NTT}. 

In this note, we are interested in the symmetrization mapping $\pi\colon \Omega_n\to \mathbb{G}_n$ and prove that the rank of the derivative $\pi'(B)$ is equal to the degree of the polynomial of $B$ for any matrix $B\in \Omega_n.$ This result explains  in part why the lifting problem is always realizable in \cite{HMY}.

\section{Content of the research}

\subsection{Notations and local lifting problem}

For $\varepsilon>0,$ denote by $$\mathbb{D}_{\varepsilon} = \{z\in \mathbb{C}~:~ |z|<\varepsilon\}$$ the disc of radius $\varepsilon.$

\begin{defn}   
For a vector $v = (v_1, v_2, \ldots, v_n)\in \mathbb{C}^n,$ put $$P_{[v]}(t) = \sum_{j=0}^n (-1)^jv_j t^{n-j}$$ where $v_0 =1$ conventionally. This notation ensures that $$P_{\pi(M)}(t) = \det(tI-M)$$ for $M\in \mathbb{C}^{n,n}.$ \end{defn}

\begin{defn}[local lifting problem] Given a holomorphic mapping $\varphi\colon \mathbb{D}_{\varepsilon}\to \mathbb{G}_n$ and $B\in \Omega_n$ such that $\varphi(0) = \pi(B).$ The \textit{local lifting problem} for $\varphi$ asks whether there exists a holomorphic mapping $\Phi\colon \mathbb{D}_{\varepsilon'}\to \Omega_n$ such that $0<\varepsilon'<\varepsilon,$ $\Phi(0) = B$ and $\pi \circ \Phi = \varphi.$
\end{defn}

\subsection{Main result} 

We will prove the following result in the remaining sections.

\begin{thm} \label{thm:main_result} Consider the symmetrization mapping $\pi\colon \Omega_n\to \mathbb{G}_n.$    Then for any $B\in \Omega_n,$ we have the rank of the derivative of $\pi$ at $B$ is equal to the degree of the minimal polynomial of $B.$
\end{thm}

We recall the definition of minimal polynomial of a square matrix $B\in \mathbb{C}^{n,n}.$

\begin{defn}
    Let $B\in \mathbb{C}^{n,n}$ be a square matrix and $P(t)$ be a complex polynomial in $t.$ We say that $P(t)$ is the minimal polynomial of $B$ if $P(t)$ is of the smallest positive degree such that $P(B) = 0.$
\end{defn}

To prove this result, we have to make use of Jordan normal form which are recalled in the next section.

\subsection{Jordan normal form and notations}

Firstly, by conjugation, we can suppose $B$ is in the Jordan normal form, i.e., $B$ is in the form of block matrix where each block in the diagonal is an elementary Jordan block corresponding to an eigenvalue. We always gather all the elementary Jordan blocks corresponding to the same eigenvalue into a big Jordan block.

It means that $B = \begin{bmatrix}
    B_1 &  & & \\ &B_2&& \\ &&\ddots& \\ &&&B_m
\end{bmatrix}$ where $B_1, B_2,\ldots, B_m$ are big Jordan blocks corresponding to distinct eigenvalues $\lambda_1, \lambda_2, \ldots, \lambda_m.$ Temporarily, we fix an eigenvalue $\lambda_k$ and its block $B_k$ and denote by $(\lambda, B)$ (i.e. without "$k$").

So $B$ here has only one eigenvalue which is $\lambda$ and $B$ is of size $m_{\lambda}$ which is the algebraic multiplicity of $\lambda.$ 

We suppose $B$ is of the form $$\begin{bmatrix}
    \lambda & b_{1,2} & 0 & & \\ &\lambda&b_{2,3}&& \\&&\ddots&& \\ &&&\lambda&b_{m_{\lambda}-1, m_{\lambda}} \\ &&&& \lambda
\end{bmatrix}$$ where the entries $b_{j-1,j}\in \{0,1\}$ and all the blank entries are equal to zero.

Then we follow \cite{NTT} and put $$F_0^{\lambda} = \{j~:~ b_{j-1,j} = 0\} = \{1=b_1 < b_2 <\ldots < b_s\}$$ and suppose $$b_2-b_1 \leq b_3 -b_2\leq \ldots \leq b_{s+1}-b_s$$ where $b_{s+1} = m_{\lambda}+1.$ This means we arrange the elementary Jordan blocks in $B$ in an increasing order of sizes from the left to the right. 

Next put $d_i^{\lambda} = 1+ \#(F^{\lambda}_0\cap [m_{\lambda-i+2..m_{\lambda}}]).$

Then if $\Phi\colon \mathbb{D}_{\varepsilon}\to \Omega_n$ is a holomorphic mapping with $\Phi(0)= B,$ and if $\varphi = \pi \circ \Phi,$ then we have \begin{equation} P^{(k)}_{[\varphi(\zeta)]}(\lambda) = O(\zeta^{d^{\lambda}_{m_{\lambda}-k}}) \mbox{ for }0\leq k\leq m_{\lambda}-1\label{eq:Trao}\end{equation} where $O$ is the big O Landau notation and $P^{(k)}_{[\varphi(\zeta)]}$ is the $k-$th derivative of the polynomial with respect to its own variable (so $\zeta$ is the variable of $\varphi$). This computation follows from the fact that $$P^{(k)}_{[\varphi(\zeta)]}(\lambda) = k!\sigma_k(\lambda I-\Phi(\zeta)).$$ 

The estimate in \eqref{eq:Trao} was first stated and proved in \cite{PT}.

\subsection{Proof of the main result}

Now we present the proof of Theorem \ref{thm:main_result}. As before, $B$ can be replaced by any similar matrix, so we can suppose $B$ is in the Jordan normal form. Denote by $m$ the degree of the minimal polynomial of $B.$ For each eigenvalue $\lambda$ (without counting multiplicities) of $B,$ denote by $s_{\lambda}$ the size of the biggest elementary Jordan block corresponding to $\lambda$ appearing in $B.$  Then we have $$m = \sum_{\lambda \in \mathrm{Sp}(B)} s_{\lambda}$$ where $Sp(B)$ is the spectrum of $B.$ For the proof of this result, cf \cite[Chapter 2, Corollary 2.3.12, page 31]{W} or for the statement of the result (without proof), cf the Vietnamese textbook \cite[Chapter 5, Proposition 5.5.4, page 139]{DDT}.

To prove the equality $\mathrm{rank}(\pi'(B)) = m,$ we prove two inequalities and divide the proof into two parts.

\subsubsection{The rank is less than or equal to the degree}Firstly, we prove $\mathrm{rank}(\pi'(B)) \leq m.$ To do it, we prove that the image of $\pi'(B)$ in $\mathbb{C}^n$ is a vector subspace of dimension $\leq m.$ We realise this by showing that this space is the solution set of a system of homogeneous linear equations with explicit coefficients as follows.

Put $\Phi(\zeta) = B + \zeta M$ for any $M\in \mathbb{C}^{n,n}$ and $\varphi = \pi \circ \Phi.$ Then we have $$\varphi'(0) = (\pi\circ \Phi)'(0) = \pi'(B)M.$$ Therefore $\varphi'(0)$ can be any vector in the image of $\pi'(B).$

By \eqref{eq:Trao}, we have $$P^{(k)}_{[\varphi(\zeta)]}(\lambda) = O(\zeta^{d^{\lambda}_{m_{\lambda}-k}}) \mbox{ for }0\leq k\leq m_{\lambda}-1.$$

Recall that $P_{[\varphi(\zeta)]}(t) = \sum_{j=0}^{n}(-1)^j\varphi_j(\zeta)t^{n-j}$ with convention $\varphi_0(\zeta) = 1.$ Put $v(t) = (-t^{n-1}, t^{n-2}, \ldots, (-1)^{n-1}t, (-1)^n).$ For two complex vectors $u = (u_1, u_2, \ldots, u_n)$ and $v = (v_1, v_2, \ldots, v_n)\in \mathbb{C}^n,$ put $$u\cdot v = \sum_{j=1}^n u_jv_j.$$ Then we have $$P^{(k)}_{[\varphi(\zeta)]}(\lambda) = v^{(k)}(\lambda)\cdot \varphi(\zeta).$$ Therefore $$\frac{d}{d\zeta}P^{(k)}_{[\varphi(\zeta)]}(\lambda) = v^{(k)}(\lambda)\cdot \varphi'(\zeta).$$ So if $d_{m_{\lambda}-k}^{\lambda}\geq 2,$ then \begin{equation}0 =\frac{d}{d\zeta}\Biggl|_{\zeta =0}P^{(k)}_{[\varphi(\zeta)]}(\lambda) = v^{(k)}(\lambda)\cdot \varphi'(0).\label{eq:nullspace}\end{equation}

Consider $d_{m_{\lambda}-k}^{\lambda}\geq 2,$ this is equivalent to $$F_0^{\lambda}\cap [m_{\lambda}-(m_{\lambda}-k)+2..m_{\lambda}] \neq \emptyset.$$ This means that $b_s\geq k+2,$ so $k\in \{0,1, \ldots, b_s-2\}.$

On the other hand, $m_{\lambda} - (b_s-1) = s_{\lambda}$ is the size of the biggest elementary Jordan block corresponding to $\lambda$ appearing in $B.$ Therefore we obtain $n-m$ homogeneous linear equations of the type \eqref{eq:nullspace} for $\varphi'(0).$

The final step of the first stage of the proof is to prove the linear independence of the coefficient vectors $$\{v^{(k)}(\lambda)~:~ \lambda\in \mathrm{Sp}(B), 0\leq k\leq m_{\lambda}-s_{\lambda}-1\}.$$

Firstly, denote by $[v_1, v_2, \ldots, v_n]$ the determinant of the matrix formed by $v_1, v_2, \ldots, v_n\in \mathbb{C}^n.$ Then we find that $$[v(\lambda_1), v(\lambda_2), \ldots, v(\lambda_n)]$$ is in fact Vandermonde's determinant, therefore $$[v(\lambda_1), v(\lambda_2), \ldots, v(\lambda_n)] = \pm \prod_{i<j}(\lambda_j-\lambda_i).$$

So $v(\lambda_1), v(\lambda_2), \ldots, v(\lambda_n)$ are linearly independent if $\lambda_1, \lambda_2, \ldots, \lambda_n$ are distinct. To deal with derivatives $v^{(k)})(\lambda),$ we make use of divided differences and we follow the notations in \cite{Boor}, precisely, $$v[\lambda_1,\lambda_2] = \frac{v(\lambda_2)-v(\lambda_1)}{\lambda_2-\lambda_1}.$$

We make first observations as follows: \begin{align*}
    [v(\lambda_1), v(\lambda_2), \ldots, v(\lambda_n)] = &  [v(\lambda_1), v(\lambda_2)- v(\lambda_1, v(\lambda_3), \ldots, v_{\lambda_n})] \\ & = (\lambda_2-\lambda_1)[v(\lambda_1), v[\lambda_1, \lambda_2], v(\lambda_3), \ldots, v(\lambda_n)].
\end{align*}

Recall that $v'(\lambda_1) = \lim_{\lambda_2\to \lambda_1} v[\lambda_1, \lambda_2].$ Therefore we get $$[v(\lambda_1), v'(\lambda_1), v(\lambda_3), \ldots, v(\lambda_n)] = \pm \prod_{\substack{i<j \\ i, j\neq 3}} (\lambda_j-\lambda_i).$$ 

By induction and combining with Genocchi-Hermite formula \cite{Boor} where $$v^{(k)}(\lambda_1) =  k! \lim_{\substack{\lambda_j\to \lambda_1\\ j\neq 1}}v[\lambda_1, \lambda_2, \ldots, \lambda_{k+1}],$$ we deduce a formula of type $$[v(\lambda_1), v'(\lambda_1), \ldots, v^{(k)}(\lambda_1), v(\lambda_2), \ldots, v(\lambda_{n-k})] = \pm 1! 2! \ldots k! \prod_{i<j}(\lambda_j-\lambda_i).$$

Therefore, we deduce that the coefficient vectors 
 $$\{v^{(k)}(\lambda)~:~ \lambda\in \mathrm{Sp}(B), 0\leq k\leq m_{\lambda}-s_{\lambda}-1\}$$ are linearly independent, and we have such $n-m$ vectors. So the solution space is of dimension $\leq m,$ it means $$\mathrm{rank}(\pi'(B))\leq m.$$

 \subsubsection{The rank is bigger than or equal to the degree} 

To realize this part of the proof, we construct $m$ linearly independent vectors in the image space of $\pi'(B)$ in the echelon form as follows. First ideas of the proof come from the last part of \cite{NTZ}.

Firstly, we suppose $B$ is in the canonical form (or rational form), i.e., $$B = \begin{pmatrix}
    C_1 & && \\ &C_2&& \\ &&\ddots& \\ &&&C_l
\end{pmatrix}$$ where $C_i's$ are all companion matrices of the invariant factors of the characteristic polynomial of $B$ and we arrange $C_i's$ in an increasing order of size, and note that the largest size is of $C_l$ and equal to $m$ the degree of the minimal polynomial. On the other hand, $C_l$ is the companion matrix of the minimal polynomial of $B,$ i.e., $C_l$ is of the form $$C_l = \begin{pmatrix}
    0 & 1 &  && \\ &0&1&& \\ &&\ddots&\ddots&\\&&&0&1 \\ -a_1&a_2&\cdots&-a_{m-1}&-a_m
\end{pmatrix}.$$

Next put $H = \begin{pmatrix}
    H_1 & && \\ &H_2&& \\ &&\ddots& \\ &&&H_l
\end{pmatrix}$ to be a diagonal block matrix of the same form as $B$ where $H_1, H_2, \ldots, H_{l-1}$ are all zero matrices and $H_l$ is of the form $$H_l = \begin{pmatrix}
    0 & 0&&&0 \\ &&\cdots&& \\ &&\cdots&& \\ &&&& \\-h_1&-h_2&\cdots&-h_{m-1}&-h_m
\end{pmatrix},$$ i.e., only the last row of $H_l$ has nonzero entries and $h_1, \ldots, h_m$ are complex variables.

Then $\sigma_k(B+\zeta H)$ is a linear polynomial of $h_m, h_{m-1}, \ldots, h_{m-k+1}$ (with nonzero coefficients on these variables) for $1\leq k\leq m.$

Therefore  if we let $h_i = 1$ and $h_j =0$ for $j\neq i$ and $1\leq j\leq m,$ the these vectors $$\pi'(B)H = \frac{d}{d\zeta}\Biggl|_{\zeta= 0} \pi(B+\zeta H)$$ form $m$ vectors in the echelon form, therefore linearly independent. It means $$\mathrm{rank}(\pi'(B)) \geq m. \qed$$

\section{Conclusion}  In this note we proved the rank of $\pi'(B)$ is equal to the degree of the minimal polynomial of $B.$ This result explains partially why we often have nice results if $B$ is a non-derogatory matrix (or cyclic matrix according to the terminology of the authors in \cite{NTT} or \cite{PT}) or in the other words, cyclic matrices are regular points of $\pi$ in the differential sense. This result partially explains why we can always lift the mapping $\varphi\colon \mathbb{D}\to \mathbb{G}_n$ locally in \cite{HMY} (although these authors prove a stronger result: they prove the global lifting).

\end{document}